\documentclass{amsart}
\newcommand{\N}{{\mathbb N}}

\newcommand{\C}{{\mathbb C}}
\newcommand{\R}{{\mathbb R}}

\newcommand{\Q}{{\mathbb Q}}
\newcommand{\Z}{{\mathbb Z}}
\newcommand{\CP}{{\mathbb CP}}

\newtheorem{theorem}{Theorem}[section]
\newtheorem{lemma}[theorem]{Lemma}

\theoremstyle{definition}
\newtheorem{definition}[theorem]{Definition}

\theoremstyle{remark}

\numberwithin{equation}{section}

\begin{document}

\title{An extension of Biran's  Lagrangian barrier theorem}

\author{Guang-Cun Lu}
\address{Department of Mathematics, Beijing Normal University, Beijing 100875,  P.R.China}
 \email{gclu@bnu.edu.cn}
\thanks{The author was supported in part by NNSF 19971045 and 10371007 of China.}


\date{July 1, 2003 and, in revised form, December 26, 2003.}

\dedicatory{{\rm (}Communicated by  Jon Wolfson{\rm )}}


\begin{abstract}
We use the Gromov-Witten invariants and a nonsqueezing theorem by
the author to affirm a conjecture by P.Biran on the Lagrangian
barriers.
\end{abstract}

\maketitle

\section{Main Results}

A K\"ahler manifold is a triple consisting of a symplectic
manifold $(M,\omega)$ and an integrable complex structure $J$
compatible with $\omega$ on $M$. If $[\omega]\in H^2(M^{2n},\Z)$
it follows from Kodaira's embedding theorem that there exists a
smooth and reduced complex hypersurface $\Sigma\subset M$ such
that its homology class $[\Sigma]\in H_{2n-2}(M)$ represents the
Poincar\'e dual $k[\omega]\in H^2(M)$ for some $k\in\N$. Following
[1] ${\mathcal P}=(M,\omega,J;\Sigma)$ is called a {\it smoothly
polarized K\"ahler manifold}. Under the conditions that either
$\dim_{\R}M\le 6$ or $\omega|_{\pi_2(M)}=0$ the following two
theorems were proved in Theorem 1.D and Theorem 4.A of [1]
respectively.\vspace{2mm}

\begin{theorem} If $(M,\omega)$ is a K\"ahler manifold with $[\omega]\in
H^2(M,\Q)$, then for every $\epsilon>0$ there exists a Lagrangian
CW-complex $\triangle_\epsilon\subset (M,\omega)$ such that every
symplectic embedding $\varphi:B(\epsilon)\to (M,\omega)$ must
satisfy $\varphi(B(\epsilon))\cap\triangle_\epsilon\ne\emptyset$.
\end{theorem}

\begin{theorem}
If ${\mathcal P}=(M,\omega, J;\Sigma)$ is a $n$-dimensional
polarized K\"ahler manifold of degree $k$ then every symplectic
embedding $\varphi: B^{2n}(\lambda):=\{
x\in\R^{2n}\,|\,|x|\le\lambda^2\}\to (M,\omega)$ with
$\lambda^2\ge\frac{1}{k\pi}$ must intersect the skeleton
$\triangle_{\mathcal P}$ associated to the polarization ${\mathcal
P}$. \end{theorem}

For the definitions of $\triangle_\epsilon$ and
$\triangle_{\mathcal P}$ the reader may refer to [1].  Actually,
such  generalizations were conjectured in Remark 4.B of [1]. As in
[1] Theorem 1.1 may be derived from Theorem 1.2.

Following [2] a Stein manifold is said to be {\it subcritical} if
it admits a plurisubharmonic function which has only critical
points of index less than half the real dimension. A polarization
${\mathcal P}=(M,\omega, J;\Sigma)$ is called {\it subcritical} if
the complement $(M\setminus\Sigma,\omega)$ is a subcritical Stein
manifold. From the proof of Theorem 1.2 we easily get the
following generalizations of Theorem F and Theorem G in
[2].\vspace{2mm}

\begin{theorem} If a closed K\"ahler
manifold $(M,\omega, J)$ admits a subcritical polarization
${\mathcal P}=(M,\omega, J;\Sigma)$ of degree $k$, then ${\mathcal
W}_G(M,\omega)\le\frac{1}{k}$ and $k\le\dim_{\C}M$. Moreover, if
the linear system of holomorphic sections of the  normal bundle to
$\Sigma$, $N_{\Sigma/M}\to\Sigma$, is base point free then
${\mathcal W}_G(M,\omega)=\frac{1}{k}\ge 1/\dim_{\C}M$. Here
${\mathcal W}_G(M,\omega)$ stands for the Gromov width of
$(M,\omega)$. \end{theorem}

\section{Proof of Theorems}

Our purpose is to prove Theorem 1.2.  The ideas are similar to
those of Biran. However, we use the theory of virtual cycles and
some techniques in [4] to compute the desired the Gromov-Witten
invariant without additional assumptions as in [1]. Then we
directly use the author's previous work on pseudo symplectic
capacities and avoid Gromov's arguments as used in [1]. For
convenience  of the reader we need to recall some related notions
in [1]. A subset $\triangle$ of a symplectic manifold $(M,\omega)$
is called an {\it embedded CW-complex} if there exists an abstract
finite and path-connected CW-complex $K$ and a homeomorphism
$i:K\to\triangle\subset M$ such that for every cell $C\subset K$
the restriction $i|_{{\rm Int}(C)}:{\rm Int}(C)\to M$ is a smooth
embedding.  When the image $i({\rm Int}(C))$ by $i$ of each cell
$C$ in $K$ is an isotropic submanifold of $(M,\omega)$ the above
embedded CW-complex $\triangle$ is called {\it isotropic}. The
dimension $\dim\triangle$ of $\triangle$ is defined as the maximum
of those of cells in $K$. An embedded isotropic CW-complex of
dimension $\frac{1}{2}\dim_{\R}M$ in $(M,\omega)$ is called a {\it
Lagrangian} CW-complex.

Another key role in [1] is the {\it standard symplectic disc
bundle}. For a closed symplectic manifold $(S,\sigma)$ with
$[\sigma]\in H^2(S;\Z)$ there exists a Hermitian line bundle
$p:L\to S$ with $c_1(L)=[\sigma]$ and a compatible connection
$\nabla$ on $L$ with curvature $R^\nabla= 2\pi i\sigma$ (cf. [5,
Prop. 8.3.1]). Denote by $\|\cdot\|$ the Hermitian metric and by
$E_L=\{v\in L\,|\,\|v\|<1\}$ the open unit disc bundle of $L$. Let
$\alpha^\nabla$ the associated {\it transgression $1$-form} on
$L\setminus 0$ with $d\alpha^\nabla=-p^\ast\sigma$, and $r$ the
radial coordinate along the fibres induced by $\|\cdot\|$. It was
shown that
$$\omega_{\rm can}:=p^\ast\sigma+ d(r^2\alpha^\nabla)$$
is a symplectic form on $E_L$ and that the symplectic type of
$(E_L, \omega_{\rm can})$ depends only on the symplectic type of
$(S,\sigma)$ and the topological type of the complex line bundle
$p:L\to S$. Moreover, $(E_L,\omega_{\rm can})$ is uniquely
characterized (up to symplectomorphism) by the following three
properties:

\begin{itemize}
\item All fibres of $p:E_L\to S$ are symplectic with respect to
$\omega_{\rm can}$ and have area $1$.
\item The restriction of
$\omega_{\rm can}$ to the zero section $S\subset E_L$ equals
$\sigma$.
 \item $\omega_{\rm can}$ is $S^1$-invariant with respect
to the obvious circle action on $E_L$.
\end{itemize}
Following [1] we call
$$p: (E_L,\omega_{\rm can})\to (S,\sigma)$$
a {\it standard symplectic disc bundle over $(S,\sigma)$ modelled
on $L$} and $p:(E_L, c\omega_{\rm can})\to (S,\sigma)$ a {\it
standard symplectic disc bundle with fibres of area $c$} for each
$c>0$.

As done in [1] this standard symplectic disc bundle can be
compactificated into a ${\CP}^1$-bundle over $(S,\sigma)$. Indeed,
let $\C$ stand for the trivial complex line bundle over $S$. Then
the direct sum $L\oplus\C$ is a complex vector bundle of rank $2$
over $S$. Denote by $p_X: X_L=P(L\oplus\C)\to S$ its projective
bundle, which is a ${\CP}^1$-bundle over $S$. It has two
distinguished sections: the zero section $Z_0=P(0\oplus\C)$ and
the section at infinity $Z_\infty=P(L\oplus 0)$. Clearly, the open
manifold $X_L\setminus Z_\infty$ is diffeomorphic to the disc
bundle $E_L$. Let $F_s$ stand for the fibre of $X_L$ at $s\in S$.
It was shown in [1] that a given $J_S\in{\mathcal J}(S,\sigma)$
and the connection $\nabla$ determine a unique almost complex
structure $J_L$ on the total space of $L$. This $J_L$ induces the
almost complex structures $J_E$ on $E$ and $J_X$ on $X_L$ again.
In particular, $J_E\in{\mathcal J}(E_L, \omega_{\rm can})$ and
$Z_0$, $Z_\infty$ and all fibres $F_s$ of $p_X: X_L\to S$ are
holomorphic with respect to $J_X$. These show that $p_X:(X_L,
J_X)\to (S,J_S)$ is an almost complex fibration with fibre
$({\CP}^1, i)$ in the following sense.\vspace{2mm}

\begin{definition}[{[1, Def.6.3.A]}]
Let $(F, J_F)$, $(B, J_B)$ and $(X, J)$ be three almost complex
manifolds. A holomorphic map $p:(X,J)\to (B, J_B)$ is called an
 {\rm almost complex fibration with fibre} $(F, J_F)$ if every
 $b\in B$ has a neighborhood $U$ and a trivialization
 $\varphi: X|_U\to U\times F$ such that for every $a\in U$ the map
 $\varphi|_{F_a}:(F_a, J|_{F_a})\to (a\times F, J_F)$ is holomorphic.
 \end{definition}

For every $0<\rho<1$ both $E_L(\rho):=\{v\in E_L\,|\,
\|v\|<\rho\}$ and $\overline{E_L(\rho)}:=\{v\in E_L\,|\,
\|v\|\le\rho\}$ are subbundles of $E_L$.\vspace{2mm}

\begin{lemma}[{[1, Lem. 5.2.A]}] Let
$p:(E_L, \omega_{\rm can})\to (S,\sigma)$ be a standard symplectic
disc bundle modeled on a Hermitian line bundle $p:L\to S$ as
above. Then there exists a diffeomorphism $f:E_L\to X_L\setminus
Z_\infty$ and a family of symplectic forms
$\{\eta_\rho\}_{0<\rho<1}$ on $X_L$ such that
\begin{enumerate}
\item $f^\ast\eta_\rho=\omega_{\rm can}$ on $\overline{E_L(\rho)}$
for every $0<\rho<1$;

\item $f$ sends the fibres of $E_L\to S$ to the fibres of
$X_L\setminus Z_\infty\to S$;

\item If $S$ is identified with the zero-section of $E_L$ then
$f(S)=Z_0$ and $p_X\circ f|_S:S\to S$ is the identity map;

\item $Z_0$, $Z_\infty$ and all fibres $F_s=p_X^{-1}(s)\, (s\in
S)$ of $p_X:X_L\to S$ are not only symplectic with respect to
$\eta_\rho$ for every $0<\rho<1$, but also holomorphic with
respect to $J_X$. Moreover, for every $0<\rho<1$ it holds that
$\eta_{\rho}|_{TZ_0}=(p^\ast_X\sigma)|_{TZ_0}$ and
$\eta_\rho|_{TZ_\infty}=c_\rho(p^\ast_X\sigma)|_{TZ_\infty}$ for
some $0<c_\rho<1-\rho^2$;

\item The area of the fibres $F_s$ satisfies
$\rho^2<\int_{F_s}\eta_\rho<1$ for every $0<\rho<1$.
\end{enumerate}
Moreover, if $J_S\in{\mathcal J}(S,\sigma)$ then $J_X\in{\mathcal
J}(X_L,\eta_\rho)$ for every $0<\rho<1$. In particular, if
$(S,\sigma,J_S)$ is K\"ahler and $L\to S$ is a holomorphic line
bundle then $(X_L,\eta_\rho, J_X)$ is K\"ahler for every
$0<\rho<1$.\end{lemma}

It was proved in Lemma 6.A of [1] that the almost complex
structure $J_X$ is regular for the class of the fibre $F\in
H_2(X_L,\Z)$ and that the space of $J_X$-holomorphic spheres in
the class $F$ is made up of exactly the fibres $F_s=p^{-1}_X(s)$,
$s\in S$. So the space ${\overline M}_{0,3}(X_L, J_X, F)$ of all
$J_X$-holomorphic stable maps in class $F$ of genus $0$ and with
$3$ marked points is nonempty. As expected we have:\vspace{2mm}

\begin{lemma} Let $F\in
H_2(X_L;\Z)$ denote the homology class of a fibre of $X_L\to S$.
Then the Gromov-Witten invariant of $(X_L,\eta_\rho)$,
$$\Psi^{(X_L,\eta_\rho)}_{F,0,3}(pt;[Z_0],[Z_\infty],pt)=1.$$
That is, $(X_L,\eta_\rho)$ is a strong $0$-symplectic uniruled
manifold in the sense of Definition 1.16 in [3].\end{lemma}

\begin{proof}Let $A=[{\CP}^1]\in
H_2({\CP^1};\Z)$ and $i:{\CP}^1\to X_L$ be the inclusion. Then
$F=i_\ast(A)$. Since the Gromov-Witten invariants are symplectic
deformation invariants, we may fix a $\eta\equiv\eta_\rho$. By
Lemma 2.2 all fibres $F_s=\pi_X^{-1}(s)$ ($s\in S$) are not only
symplectic with respect to $\eta$ but also holomorphic with
respect to $J_X$. Thus for each $s\in S$ and $x\in F_s$ the
tangent space $T_xF_s$ is a symplectic and $J_X(x)$-invariant
subspace of $(T_xX_L,\eta_x)$. Let $H_x$ denote the
$\eta_x$-orthogonal complement of $T_xF_s$. Then $T_xX_L=H_x\oplus
T_xF_s$ and $Dp_X(x)|_{H_x}:H_x\to T_sS$ is an isomorphism. Note
that $J_X$ is also compatible with $\eta$. We have that
$\eta(J_X\xi, J_X\eta)=\eta(\xi,\eta)$ for all $\xi,\eta\in TX_L$.
This implies that $H_x$ is also a $J_X(x)$-invariant subspace of
$T_xX_L$.  Consequently, $J_X(x)$ preserve the splitting
$T_xX_L=H_x\oplus T_xF_s$. Moreover, the projection $p_X: (X_L,
J_X)\to (S, J_S)$ is holomorphic(cf. [1, \S6.3]). Hence the almost
complex structure $J_X$ on $X_L$ is {\it fibered} in the sense of
Definition 2.2 in [4].

 Since any stable map has connected image set it
is easily checked that each stable map in
 ${\overline M}_{0,3}(X_L, J_X, F)$ must also entirely lie in
 a fibre of $X_L$. It follows from this and Lemma 4.3 in [4] that
for any representative $\tilde\tau=(\Sigma, \tilde h)$ of $\tau\in
{\overline M}_{0,3}(X_L, J_X, F)$ the cokernel of
$D\bar\partial_{J_X}(\tilde h)$ can be spanned by elements of the
space
$${\mathcal L}_{\tilde h}^V=L^{1,p}\Bigr(\wedge^{0,1}
({\tilde h}^\ast(T_{vert}X_L))\Bigl),$$
where $V=T_{vert}X_L$ is the vertical tangent bundle of $TX_L$. By
Proposition 4.3 in [4] one can choose $R$ and the embeddings $e$
in the construction of the virtual moduli cycle ${\overline
M}^\nu_{0,3}(X_L, J_X, F)$ of ${\overline M}_{0,3}(X_L, J_X, F)$
such that $e_{\tilde\tau}(\nu)\in{\mathcal L}^V_{\tilde h}$  for
all $\nu\in R$ and all $\tilde\tau$. That is, we can choose a
fibered pair $(J_X,\nu)$ on $X_L$ in the sense of Definition 4.4
in [4]. As proved in Proposition 4.4 of [4] such a virtual moduli
cycle ${\overline M}^\nu_{0,3}(X_L, J_X, F)$ has the following
property. For each element of it, a parameterized stable map
$(\Sigma,\tilde h)$, each component $\tilde h_k$ of $\tilde h$
satisfies an equation of the form $\bar\partial_{J_X}\tilde
h_k=\nu_k$ for some section $\nu_k$ of $\wedge^{0,1}(\tilde
h_k^\ast(T_{vert}X_L))$. This implies that each $p_X\circ\tilde
h_k:S^2\to S$ is $(i, J_S)$-holomorphic and therefore that each
element in ${\overline M}^\nu_{0,3}(X_L, J_X, F)$ has the image
contained in a single fibre of $X_L$. For any $s\in S$ we identify
$(F_s, J_X|_{F_s})$ with $(\CP^1, i)$. Then the subset ${\overline
M}^\nu_{0,3}(X_L, J_X, F)_s$ of ${\overline M}^\nu_{0,3}(X_L, J_X,
F)$ consisting of stable maps with image in $F_s$ regularizes the
space ${\overline M}_{0,3}(\CP^1, i, A)$ in the sense of
Definition 4.3 in [4]. Hence
$$\Psi^{(X_L,\eta)}_{F,0,3}(pt;[Z_0],[Z_\infty],pt)=
\Psi^{(\CP^1,\omega_{FS})}_{A,0,3}(pt;pt, pt, pt)=1.$$ Here we
have used $pt$ to stand for the point classes in the different
spaces. Lemma 2.3 is proved.\end{proof}

Recall that for a closed symplectic manifold $(M,\omega)$ and
homology classes $\alpha_0,\alpha_\infty\in H_\ast(M,\Q)$ we in
[3, Def.1.8] defined a number
$${\rm GW}_g(M,\omega;\alpha_0,\alpha_\infty)\in (0, +\infty]$$
by the infimum of the $\omega$-areas $\omega(A)$ of the homology
classes $A\in H_2(M;\Z)$ for which the Gromov-Witten invariant
$\Psi_{A, g,
m+2}(C;\alpha_0,\alpha_\infty,\beta_1,\cdots,\beta_m)\ne 0$ for
some homology classes $\beta_1,\cdots,\beta_m\in H_\ast(M;\Q)$ and
$C\in  H_\ast(\overline{\mathcal M}_{g, m+2};\Q)$ and integer
$m>0$. Moreover, in Definition 1.25 of [3] we also defined another
number ${\rm GW}(M,\omega)\in (0, +\infty]$ by
$${\rm GW}(M,\omega)=\inf{\rm GW}_g(M,\omega; pt,\alpha),$$
where the infimum is taken over all nonnegative integers $g$ and
all homology classes $\alpha\in H_\ast(M;\Q)\setminus\{0\}$ of
degree $\deg\alpha\le\dim M-1$. Theorem 1.26 in [3] claimed that
$${\mathcal W}_G(M,\omega)\le{\rm GW}(M,\omega)\leqno(1)$$
for any symplectic uniruled manifold $(M,\omega)$. \vspace{2mm}

\noindent{\it Proof of Theorem 1.2}.\hspace{2mm} By the
decomposition Theorem 2.6.A in [1],
$(M\setminus\triangle_{\mathcal P},\omega)$ is symplectomorphic to
the standard symplectic disc bundle
$(E_{N_\Sigma},\frac{1}{k}\omega_{\rm can})\to
(\Sigma,\sigma=k\omega|_{\Sigma})$ over $\Sigma$ which is modeled
on the normal bundle $N_\Sigma$ and with fibres of area $1/k$.
Assume that $\varphi(B(\lambda))\cap\triangle_{\mathcal
P}=\emptyset$ for some symplectic embedding $\varphi:
B^{2n}(\lambda)\to (M,\omega)$ with $\lambda^2\ge\frac{1}{k\pi}$.
Then there exists a symplectic embedding $\psi: B(\lambda)\to
(E_{N_\Sigma},\frac{1}{k}\omega_{\rm can})$. Since
$\psi(B(\lambda))$ is compact there exists a positive number
$\rho\in (0, 1)$ such that $\psi(B(\lambda))\subset
E_{N_\Sigma}(\rho)$.  Lemma 2.2  gives a symplectic embedding from
$\psi(B(\lambda))$ into $(X_{N_\Sigma}, \frac{1}{k}\eta_\rho)$.
Hence
$${\mathcal W}_G(X_{N_\Sigma},
\frac{1}{k}\eta_\rho)\ge\pi\lambda^2\ge 1/k.\leqno(2)$$

On the another hand, by  the definition of ${\rm GW}(M,\omega)$
and Lemmas 2.2, 2.3 we have
$${\rm GW}(X_{N_\Sigma},\eta_\rho)\le \int_{F_s}\eta_\rho<1$$
for this  $\rho$. This and (1) together give
$${\mathcal W}_G(X_{N_\Sigma},
\frac{1}{k}\eta_\rho) \le{\rm
GW}(X_{N_\Sigma},\frac{1}{k}\eta_\rho)\le
\frac{1}{k}\int_{F_s}\eta_\rho<\frac{1}{k},$$
 which contradicts to (2). Theorem 1.2 is proved.
 \hfill$\Box$\vspace{2mm}

\noindent{\bf Acknowledgements}.\hspace{2mm}I am grateful to Kai
Cieliebak for sending me some preprints. I would like to thank the
referees for several advice for improving the presentation,
pointing out a careless mistake in the statement of Theorem 1.3 in
the previous version of the paper.

\end{document}